\documentclass[12pt]{amsart}

\usepackage{amssymb}
\usepackage{tikz}   

\usepackage{enumerate}

\usepackage{graphicx}

\makeatletter
\@namedef{subjclassname@2020}{%
  \textup{2020} Mathematics Subject Classification}
\makeatother

\newtheorem{thm}{Theorem}[section]

\newtheorem{prop}[thm]{Proposition}

\theoremstyle{definition}

\newtheorem{rem}[thm]{Remark}
\newtheorem{exa}[thm]{Example}

\numberwithin{equation}{section}

\begin{document}


\baselineskip=17pt



\title[A note on Carleson sets and Carleson windows]{A note on Carleson sets and Carleson windows}

\author[Mohammad Ansari]{Mohammad Ansari}
\dedicatory{}


\date{}
\begin{abstract} Let $\Bbb D$ be the open unit disk in the complex plane. 
In this note, for the Carleson set $S(b,h)$ and the Carleson 
window $W(b,h)$ which are particular subsets of $\Bbb D$, we find conditions on $h$ under which we 
have $W(b,h/c)\subset S(b,h)\subset W(b,ch)$ for some $c>1$. 
\end{abstract}


\maketitle

\section{\bf Introduction}
In this note, we give a result (Theorem $2.1$) which shows the possible and impossible inclusions between Carleson sets and Carleson windows which are particular subsets
of the open unit disk $\Bbb D=\{z\in \Bbb C: |z|< 1\}$ in the complex plane $\Bbb C$.
\par Consider the open unit disk $\Bbb D$ in the complex plane (the disk centered at $O$ with radius $1=OM=ON$ in Figure $1$),
choose a point $b$ on the boundary of $\Bbb D$, i.e., on the unit circle $\Bbb T$, and let $0<h<1$. 
\vspace*{.5cm}
\begin{center}
\begin{tikzpicture}
\draw[thick](1.5,0) circle [radius=1.5];
\draw[thick](0,0) circle [radius=1];
\draw[thick](1.5,0) circle [radius=.5];
 \node at (-1.15,0) (nodeA){} ;
    \node at (.1,0) (nodeB) {$\;b$};
    \node at (-.1,0) (nodeC) {};
   
    \node at (1.85,0) (nodeO){$O$};
    
    \draw[] (nodeA) -- (nodeB)--(nodeC)--(nodeO);
    \node at (.01,1.15) (nodeM){$M$};
    \node at (1.8,-.1) (nodeO){};
    \draw[] (nodeM)--(nodeO);
    \node at (.01,-1.15) (nodeN){$N$};
    \node at (1.8,.1) (nodeR){};
    \draw[] (nodeN)--(nodeR);
     \node at (1.05,.6) (nodeP){\;\;$P$};
     \node at (1.05,-.7)(nodeQ){\;\;$Q$};
\end{tikzpicture}
\end{center}
$$\textnormal{\small{Fig. 1}}$$
\par In the paragraph before 
Theorem $2.33$ and the paragraph after Definition $2.34$ of \cite{cm}, the {\it Carleson set} 
and the {\it Carleson window} (with parameters $b, h$) have been introduced (respectively) by $$S(b,h)=\{z\in \Bbb D: |z-b|<h\},$$ and 
$$W(b,h)=\{z\in \Bbb D: 1-h<|z|<1 \;\textnormal{and}\; z/|z|\in S(b,h)\}.$$
But, since no unimodular number $z/|z|$ could be in $S(b,h)\subset \Bbb D$, it seems that the definition of $W(b,h)$ must be modified to
$$W(b,h)=\{z\in \Bbb D: |z|>1-h \;\textnormal{and}\; z/|z|\in \overline{S(b,h)}\}.$$
\par Let $\Bbb D_1$ and $\Bbb D_2$ be, respectively, the open disk with radius $h$ centered at $b$, and the open disk with radius $1-h$ ($=OP=OQ$) centered at $O$. Meanwhile,
suppose that the circles $\Bbb T_1$ and $\Bbb T_2$ are, respectively, the boundaries of $\Bbb D_1$ and $\Bbb D_2$. Then it is clear that $S(b,h)=\Bbb D_1\cap \Bbb D$ 
and $W(b,h)$ is the subset of $\Bbb D$ whose boundary consists of the arc $MN$ of $\Bbb T$, the arc $PQ$ of $\Bbb T_2$, and the line segments $MP$ and $QN$.
\par  In the following proposition, we see that the inclusion $S(b,h)\subset W(b,h)$ cannot hold. 
 Note that if $A, B$ are the endpoints of the line segment $AB$, we also write $AB$ to refer to the length of $AB$ in our
computations. Meanwhile, the only lowercase letter which is used here is $b\in \Bbb T$. So, for example, $bM$ may refer to both the line segment $bM$
and its length. 
\begin{prop} Assume that $b\in \Bbb T$ and $0<h<1$. Then the inclusion $S(b,h)\subset W(b,h)$ cannot hold.
\begin{proof} We show that $OM$ and $ON$ are not tangent to $\Bbb T_1$ (see Figure $1$). Indeed if they were, regarding the power
of the point $O$ with respect to the circle $\Bbb T_1$, we could write
 $$1=OM^2=ON^2=1-h^2<1,$$ which is a contradiction. Therefore, the half-lines $\vec{OM}$ and 
$\vec{ON}$ must intersect
$\Bbb T_1$ at other points different from $M$ and $N$. If we name them $M'$ and $N'$, then it is clear that $M', N'\in \Bbb D$. In fact, $OM'=OM'.OM=1-h^2<1$ and similarly $ON'<1$. 
Now, it is clear that $S(b,h)$ cannot be a subset of $W(b,h)$.
\end{proof}
\end{prop}
The above result motivated us to search for those constants $c>1$ for which $S(b,h)\subset W(b,ch)$. As the authors in \cite{cm} have mentioned, we may sometimes need to 
replace $S(b,h)$ with $W(b,h)$ or vice versa in some relevant computations. Then we may need to have 
the inclusions $W(b,h/c)\subset S(b,h)$ and $S(b,h)\subset W(b,ch)$ for some $c>1$. In Section $2$, we investigate the possibilities of these inclusions.
\section{\bf The inclusions between Carleson sets and windows}
In the next theorem, we discuss the possible and impossible inclusions between $W(b,h/c)$, $S(b,h)$, and $W(b,ch)$ for $c>1$.
\begin{thm} For any fixed $b\in \Bbb T$ the following hold.\\
\hspace*{.3cm} \textnormal{(i)} For any $0<h<1$ there is a ray $R_h\subset (1,\infty)$ such that for any $c\in R_h$ we have $W(b,h/c)\subset S(b,h)$.\\
\hspace*{.2cm} \textnormal{(ii)} If $0<h$, then there is some $c>1$ such that $S(b,h)\subset W(b,ch)$ if and only if $h<\sqrt{3}/2$. In that case, the desired constants
$c$ make a bounded interval $I_h$.\\
\hspace*{.1cm} \textnormal{(iii)} If $0<h$, then there is some $c>1$ such that $$W(b,h/c)\subset S(b,h)\subset W(b,ch)$$ if and only if $h<\sqrt{3}/2$. In that case, the desired 
constants $c$ make a bounded interval $I_h$.
\begin{proof} \textnormal{(i)}. In view of Figure $2$, assume that for some $c>1$, the disks centered at $b$ have radii
$h/c$ and $h$ and we have $W(b,h/c)\subset S(b,h)$.
\vspace*{.5cm}
\begin{center}
\begin{tikzpicture}
\draw[thick](1.5,0) circle [radius=1.5];
\draw[thick](0,0) circle [radius=1];
\draw[thick](0,0) circle [radius=1.4];
\draw[thick](1.5,0) circle [radius=.5];
 \node at (-1.15,0) (nodeA){} ;
    \node at (.1,0) (nodeB) {$\;b$};
    \node at (-.1,0) (nodeC) {};
   
    \node at (1.85,0) (nodeO){$O$};
    
    \draw[] (nodeA) -- (nodeB)--(nodeC)--(nodeO);
    \node at (.01,1.15) (nodeM){$M$};
    \node at (1.8,-.1) (nodeO){};
    \draw[] (nodeM)--(nodeO);
    \node at (.01,-1.15) (nodeN){$N$};
    \node at (1.8,.1) (nodeR){};
    \draw[] (nodeN)--(nodeR);
     \node at (1.05,.6) (nodeP){\;\;$P$};
     \node at (1.05,-.7)(nodeQ){\;\;$Q$};
\end{tikzpicture}
\end{center}
$$\textnormal{\small{Fig. 2}}$$
It is clear that the constant $c>1$ must be large enough such that the points $P, Q$ of the boundary of $W(b, h/c)$
satisfy $bP\le h$ and $bQ\le h$. 
So, regarding the law of cosines in the triangles $bOM$ and $bOP$, we must have $$(h/c)^2=bM^2=2-2\cos\theta$$ and $$h^2\ge bP^2=1+(1-h/c)^2-2(1-h/c)\cos\theta,$$ where $\theta$
is the angle between $OM$ and $Ob$ (remember that $OP=1-h/c$).
Then we obtain $2c-c^3\le h$. Now, it is clear that for any $0<h<1$ there is some $c>1$ such that $2c-c^3\le h$. Moreover, 
if $f(x)=2x-x^3$ ($x>1$), then, for any given $h\in (0,1)$, the constants $c$ for which $W(b,h/c)\subset S(b,h)$ make the ray $R_h=[f^{-1}(h), \infty)$ (the continuous 
 function $f$ is strictly decreasing and $\lim_{x\to \infty}f(x)=-\infty$). 
\par \textnormal{(ii)}. If for $h>0$ and $c>1$ we have that $S(b,h)\subset W(b,ch)$, then, regarding Figure $3$ in which the disks centered at $b$ have radii $h$ and $ch$, we must have 
$h\le bR$. Here $R$ is the midpoint of the chord $M'M$ and so $bR$ is perpendicular to $M'M$ (to avoid crowdedness in the figure, $R$ and $bR$ are not shown).
 If $M'M=2x$, then we have $h^2\le bR^2=bM^2-x^2=(ch)^2-x^2$. On the other hand,
$$x=(1-OM')/2=(1-OM'.OM)/2 =[1-(1-(ch)^2)]/2=c^2h^2/2.$$ Hence, $h^2\le c^2h^2-c^4h^4/4$. If we put $y=c^2$, we must have $h^2y^2-4y+4\le 0$.
If we solve this inequality, we will have $$(\sqrt{2}/h)\sqrt{1-\sqrt{1-h^2}}\le c\le (\sqrt{2}/h)\sqrt{1+\sqrt{1-h^2}}.$$
But, we should have $ch<1$ because $ch$ is the parameter of the Carleson window $W(b, ch)$. Thus, $c<1/h$ and then we have to replace the above upper bound by the smaller one $1/h$.
On the other hand, it is easy to see that, for all $h\in (0,1)$, we have  $1<(\sqrt{2}/h)\sqrt{1-\sqrt{1-h^2}}$. 
Therefore,  the desired constants $c$ should be selected from the interval $I_h=\big{[}(\sqrt{2}/h)\sqrt{1-\sqrt{1-h^2}}, 1/h\big{)}$, and it can be readily verified that 
$I_h\neq \emptyset$ if and only if $h<\sqrt{3}/2$.
\vspace*{.5cm}
\begin{center}
\begin{tikzpicture}
\draw[thick](1.5,0) circle [radius=1.5];
\draw[thick](0,0) circle [radius=1];
\draw[thick](0,0) circle [radius=.7];
\draw[thick](1.5,0) circle [radius=.5];
 \node at (-1.15,0) (nodeA){} ;
    \node at (.1,0) (nodeB) {$\;b$};
    \node at (-.1,0) (nodeC) {};
   
    \node at (1.85,0) (nodeO){$O$};
    
    \draw[] (nodeA) -- (nodeB)--(nodeC)--(nodeO);
    \node at (.01,1.15) (nodeM){$M$};
    \node at (1.8,-.1) (nodeO){};
     \node at (1.1,.8) (nodeM'){$M'$};
        \node at (.8,.58) (nodeL){};
        \node at (.95,.48) (nodeS){};
    \draw[] (nodeM)--(nodeS)--(nodeL)--(nodeO);
    \node at (.01,-1.15) (nodeN){$N$};
    \node at (1.8,.1) (nodeR){};
    \draw[] (nodeN)--(nodeR);
     \node at (1.05,.6) (nodeP){};
     \node at (1.05,-.7)(nodeQ){};
\end{tikzpicture}
\end{center}
$$\textnormal{\small{Fig. 3}}$$
\par \textnormal{(iii)}. To have a constant $c>1$ satisfying $W(b,h/c)\subset S(b,h)\subset W(b,ch)$, we need to apply all obtained restrictions
from \textnormal{(i)} and \textnormal{(ii)}.
So, we should have $2c-c^3\le h$ and $(\sqrt{2}/h)\sqrt{1-\sqrt{1-h^2}}\le c<1/h$ simultaneously. By \textnormal{(ii)} the latter restriction holds if and only if $h<\sqrt{3}/2$.
Now, the inequality $c<1/h$ is equivalent to $2/h-1/h^3<2c-c^3$
and hence the mentioned restrictions have common solutions if and 
only if $f(1/h)=2/h-1/h^3<h$ (where $f$ is the function introduced in the proof of \textnormal{(i)}), which is true for all $0<h<\sqrt{3}/2$. This shows that
$f^{-1}(h)<1/h$. To find the interval of the desired constants $c$, note that
$2c-c^3\le h$ implies $f^{-1}(h)\le c$. Thus, the desired constants $c$ make the interval $I_h=[g(h), 1/h)$ where $g(h)=\max\{f^{-1}(h),(\sqrt{2}/h)\sqrt{1-\sqrt{1-h^2}}\}$.
\end{proof}
\end{thm}
\begin{rem} Let us have a note on the function $g$. It is easily seen that the function $k(h)=(\sqrt{2}/h)\sqrt{1-\sqrt{1-h^2}}$ is strictly increasing on $(0,\sqrt{3}/2)$.
Moreover, we know that $k(h)>1$ for all $h\in (0,\sqrt{3}/2)$ (see the proof of \textnormal({ii})) and $f$ is strictly decreasing on $(1,\infty)$.
Thus, the function $F(h)=f(k(h))-h$ is also strictly decreasing on $(0,\sqrt{3}/2)$. On the other hand, one can readily verify that
$F(0.82)>0$ and $F(0.83)<0$, and so $F$ has a unique root $h_0$ in $(0,\sqrt{3}/2)$ such that $0.82<h_0<0.83$
 (a computer aided computation shows that $h_0\approx 0.82056$). Therefore, $F(h)>0$ for $0<h<h_0$ and $F(h)\le 0$ for $h_0\le h<\sqrt{3}/2$, 
 and so $g(h)=f^{-1}(h)$ for $0<h<h_0$ and
$g(h)=k(h)$ for $h_0\le h<\sqrt{3}/2$. 
\end{rem}
\begin{exa} If $h=0.5$, then $h<h_0$ and so $I_h=[f^{-1}(h), 1/h)$. But $f^{-1}(0.5)<1.27$ and hence we have $I_h\supset [1.27,2)$. 
Then, for any $c\in [1.27,2)$ we have $W(b,0.5/c)\subset S(b, 0.5)\subset W(b, 0.5c)$. As another sample, if 
$h=0.85$, then $h>h_0$ and hence $I_h=[k(h),1/h)\supset [1.15, 1.17]$. So, if $1.15\le c\le 1.17$, then we have $W(b,0.85/c)\subset S(b, 0.85)\subset W(b,0.85c)$.
\end{exa}
 
\vspace*{.5cm}

\end{document}